\documentclass{amsart}
\usepackage{amsmath,amsthm,graphicx}
\usepackage[all]{xy}

\theoremstyle{plain}
\newtheorem{theorem}{Theorem}
\newtheorem{lemma}{Lemma}
\newtheorem{proposition}{Proposition}
\newtheorem{corollary}{Corollary}

\theoremstyle{definition}
\newtheorem{definition}{Definition}
\newtheorem{example}{Example}

\newcommand\Z{{\mathbb{Z}}}

\begin{document}

\title{Frobenius {P}roblem and dead ends in integers}

\author{Zoran \v{S}uni\'c}
\thanks{Partially supported by NSF grant DMS-0600975}
\address{Department of Mathematics, Texas A\&M University, College Station, TX 77843-3368, USA}



\begin{abstract}
Let $a$ and $b$ be positive and relatively prime integers. We show
that the following are equivalent: (i) $d$ is a dead end in the
(symmetric) Cayley graph of $\Z$ with respect to $a$ and $b$, (ii)
$d$ is a Frobenius value with respect to $a$ and $b$ (it cannot be
written as a non-negative or non-positive integer linear combination
of $a$ and $b$), and $d$ is maximal (in the Cayley graph) with
respect to this property. In addition, for given integers $a$ and
$b$, we explicitly describe all such elements in $\Z$. Finally, we
show that $\Z$ has only finitely many dead ends with respect to any
finite symmetric generating set. In the appendix we show that every
finitely generated group has a generating set with respect to which
dead ends exist.
\end{abstract}

\maketitle


\section{Introduction}

We first describe the variant of Frobenius Problem that is in our
interest.

\begin{definition}[Frobenius values]
Let $S$ be a set of positive integers whose greatest common divisor
is $1$. An integer $n$ is termed \emph{positively generated} with
respect to $S$ if it is a non-negative integer linear combination of
the elements in $S$, \emph{negatively generated} if it is a
non-positive integer linear combination of the elements in $S$, and
is termed \emph{Frobenuis value} (with respect to $S$) otherwise.
\end{definition}

It is known that, for any set $S$ of positive integers with greatest
common divisor $1$, there exist only finitely many Frobenius values.
Frobenius Problem (also known as Linear Diophantine Problem of
Frobenius) asks to find the largest Frobenius value for a given $S$.
The largest Frobenius value is called the \emph{Frobenius number} of
$S$. It is known that, for $S=\{a,b\}$, where $a$ and $b$ are
positive and relatively prime integers, the Frobenius number is
$ab-a-b$. No explicit formula exists when $S$ consists of at least
three distinct numbers. On the positive side, upper bounds exist
(see~\cite{beck-d-r:frobenius} and~\cite{fukshansky-r:frobenius} for
some estimates and further references) as do polynomial time
algorithms determining the Frobenius number for sets $S$ of fixed
size~\cite{kannan:frobenius}.

While Frobenius Problem has a long history, the notion of a dead end
is fairly recent. It appears explicitly in the work of
Bogopol$'$ski{\u\i}~\cite{bogopolskii:bi-lipschitz}, who shows that,
for a given non-elementary hyperbolic group with a given finite
generating set, there exists a uniform bound on the depth of the
dead ends in the group. Various results regarding dead ends in
Thompson's group $F$, lamplighter groups, solvable groups, finitely
presented groups, residually finite groups, etc., appear in the
works of Cleary, Guba, Riley, Taback and
Warshall~\cite{cleary-t:thompson-cayley,cleary-t:lamplighter-dead,guba:solvable-dead,cleary-r:fp-dead,riley-w:notinvariant,warshall:many,warshall:lattices}.

\begin{definition}[Word length]
Let $G$ be a group generated by a finite set $S$. For an element $g$
in $G$, define the \emph{word length} (or simply \emph{length}) of
$g$ with respect to $S$, denoted by $\ell_S(g)$ (or simply by
$\ell(g)$ when $S$ is assumed), to be the shortest length of a group
word over $S$ representing $g$, i.e.,
\[
  \ell_S(g) = \min\{\ k \mid g=s_1s_2\dots s_k, \text{ for some }
                   s_1,\dots,s_k \in S \cup S^{-1} \ \}.
\]
\end{definition}

\begin{definition}[Cayley graph]
Let $G$ be a group generated by a finite set $S$. The (symmetric)
\emph{Cayley graph} of $G$ with respect to $S$ is the graph
$\Gamma(G,S)$ whose vertices are the elements of $G$ and in which
two vertices $g$ and $h$ are connected by an edge if and only if $g
= hs$, for some $s$ in $S\cup S^{-1}$.
\end{definition}

\begin{example}
The Cayley graph of $\Z$ with respect to $S=\{3,5\}$ is given in
Figure~\ref{f:z35}.
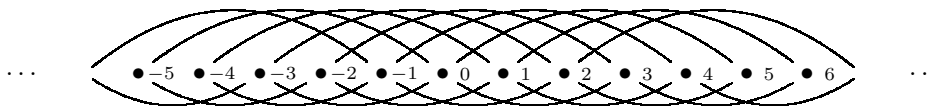
\begin{figure}[!ht]
\[
\xymatrix@C=12pt@R=10pt{
 \\
 \dots &
 \ar@{-}@/^2pc/[rrrrr] \ar@{-}@/_1pc/[rrr]&
 \bullet \ar@{}[r]|<<{-5} \ar@{-}@/^2pc/[rrrrr] \ar@{-}@/_1pc/[rrr]&
 \bullet \ar@{}[r]|<<{-4} \ar@{-}@/^2pc/[rrrrr] \ar@{-}@/_1pc/[rrr]&
 \bullet \ar@{}[r]|<<{-3} \ar@{-}@/^2pc/[rrrrr] \ar@{-}@/_1pc/[rrr]&
 \bullet \ar@{}[r]|<<{-2} \ar@{-}@/^2pc/[rrrrr] \ar@{-}@/_1pc/[rrr]&
 \bullet \ar@{}[r]|<<{-1} \ar@{-}@/^2pc/[rrrrr] \ar@{-}@/_1pc/[rrr]&
 \bullet \ar@{}[r]|<<{0}  \ar@{-}@/^2pc/[rrrrr] \ar@{-}@/_1pc/[rrr]&
 \bullet \ar@{}[r]|<<{1}  \ar@{-}@/^2pc/[rrrrr] \ar@{-}@/_1pc/[rrr]&
 \bullet \ar@{}[r]|<<{2}  \ar@{-}@/^2pc/[rrrrr] \ar@{-}@/_1pc/[rrr]&
 \bullet \ar@{}[r]|<<{3}  \ar@{-}@/_1pc/[rrr]&
 \bullet \ar@{}[r]|<<{4}  \ar@{-}@/_1pc/[rrr]&
 \bullet \ar@{}[r]|<<{5} &
 \bullet \ar@{}[r]|<<{6} &
 & \dots
 \\
 &
}
\]
 \caption{The Cayley graph of $\Z$ with respect to $S=\{3,5\}$}
 \label{f:z35}
\end{figure}

\end{example}

It is clear that, for an element $g$ in $G$, the combinatorial
distance $d_S(1,g)$ between 1 and $g$ in the Cayley graph
$\Gamma(G,S)$ is precisely the word length $\ell_S(g)$ of $g$ with
respect to $S$. The elements of length $n$ in $G$ with respect to
$S$ are precisely the elements on the \emph{sphere} $\Sigma(n)$ of
radius $n$ in $\Gamma(G,S)$.

\begin{definition}[Dead end]
Let $G$ be a group generated by a finite set $S$. A \emph{dead end}
in $G$ with respect to $S$ is an element $d$ in $G$ such that
\[ \ell(ds) \leq \ell(d) \]
for every $s$ in $S \cup S^{-1}$. A \emph{strict dead end} is an
element $d$ such that $\ell(ds)<\ell(d)$, for every $s$ in $S \cup
S^{-1}$.
\end{definition}

In the Cayley graph, a dead end is a vertex $d$ for which a geodesic
path from $1$ to $d$ cannot be further extended to a longer
geodesic. In other words, a dead end of length $n$ is a vertex in
the sphere $\Sigma(n)$ of radius $n$ from which the sphere
$\Sigma(n+1)$ of radius $n+1$ cannot be reached in one step.

\begin{example}
The Cayley graph $\Gamma(\Z,S)$, for $S=\{3,5\}$, is presented in
Figure~\ref{f:z-dead} in such a way that the length of each element
is apparent from the figure (it corresponds to the level at which it
is drawn).
\begin{figure}[!ht]
\[
\xymatrix@C=14pt{
  \bullet \ar@{}[r]|<<<{-25} \ar@{-}[d]&
  \bullet \ar@{}[r]|<<<{-23} \ar@{-}[d] \ar@{-}[ld]&
  \bullet \ar@{}[r]|<<<{-21} \ar@{-}[d] \ar@{-}[ld]&
  \bullet \ar@{}[r]|<<<{-19} \ar@{-}[d] \ar@{-}[ld]&
  \bullet \ar@{}[r]|<<<{-17} \ar@{-}[d] \ar@{-}[ld] \ar@{-}[lllld]&
  &&&
  \bullet \ar@{}[l]|<<{17} \ar@{-}[d] \ar@{-}[rd] \ar@{-}[rrrrd]&
  \bullet \ar@{}[l]|<<{19} \ar@{-}[d] \ar@{-}[rd]&
  \bullet \ar@{}[l]|<<{21} \ar@{-}[d] \ar@{-}[rd]&
  \bullet \ar@{}[l]|<<{23} \ar@{-}[d] \ar@{-}[rd]&
  \bullet \ar@{}[l]|<<{25} \ar@{-}[d]
  \\
  \bullet \ar@{}[r]|<<<{-20} \ar@{-}[d]&
  \bullet \ar@{}[r]|<<<{-18} \ar@{-}[d] \ar@{-}[ld]&
  \bullet \ar@{}[r]|<<<{-16} \ar@{-}[d] \ar@{-}[ld]&
  \bullet \ar@{}[r]|<<<{-14} \ar@{-}[d] \ar@{-}[ld]&
  \bullet \ar@{}[r]|<<<{-12} \ar@{-}[d] \ar@{-}[ld] \ar@{-}[lllld]&
  \diamond \ar@{}[r]|<<<{-4} \ar@{-}[lld] \ar@{-}[ld] \ar@{-}[d] \ar@{-}[rrd]&
  &
  \diamond \ar@{}[l]|<<{4} \ar@{-}[lld] \ar@{-}[d] \ar@{-}[rd] \ar@{-}[rrd]&
  \bullet \ar@{}[l]|<<{12} \ar@{-}[d] \ar@{-}[rd] \ar@{-}[rrrrd]&
  \bullet \ar@{}[l]|<<{14} \ar@{-}[d] \ar@{-}[rd]&
  \bullet \ar@{}[l]|<<{16} \ar@{-}[d] \ar@{-}[rd]&
  \bullet \ar@{}[l]|<<{18} \ar@{-}[d] \ar@{-}[rd]&
  \bullet \ar@{}[l]|<<{20} \ar@{-}[d]
  \\
  \bullet \ar@{}[r]|<<<{-15} \ar@{-}[rrd]&
  \bullet \ar@{}[r]|<<<{-13} \ar@{-}[rd] \ar@{-}[rrd]&
  \bullet \ar@{}[r]|<<<{-11} \ar@{-}[rd] \ar@{-}[rrd]&
  \bullet \ar@{}[r]|<<<{-9} \ar@{-}[rd]&
  \diamond \ar@{}[r]|<<<{-7} \ar@{-}[lld] \ar@{-}[rd]&
  \diamond \ar@{}[r]|<<<{-1} \ar@{-}[ld] \ar@{-}[rrd]&
  &
  \diamond \ar@{}[l]|<<{1} \ar@{-}[lld] \ar@{-}[rd]&
  \diamond \ar@{}[l]|<<{7} \ar@{-}[ld] \ar@{-}[rrd]&
  \bullet \ar@{}[l]|<<{9} \ar@{-}[ld]&
  \bullet \ar@{}[l]|<<{11} \ar@{-}[ld] \ar@{-}[lld]&
  \bullet \ar@{}[l]|<<{13} \ar@{-}[ld] \ar@{-}[lld]&
  \bullet \ar@{}[l]|<<{15} \ar@{-}[lld]
  \\
  &&
  \bullet \ar@{}[r]|<<<{-10} \ar@{-}[rrd]&
  \bullet \ar@{}[r]|<<<{-8} \ar@{-}[rd] \ar@{-}[rrd]&
  \bullet \ar@{}[r]|<<<{-6} \ar@{-}[rd]&
  \diamond \ar@{}[r]|<<<{-2} \ar@{-}[ld] \ar@{-}[rrd]&
  &
  \diamond \ar@{}[l]|<<{2} \ar@{-}[lld] \ar@{-}[rd]&
  \bullet \ar@{}[l]|<<{6} \ar@{-}[ld]&
  \bullet \ar@{}[l]|<<{8} \ar@{-}[lld] \ar@{-}[ld]&
  \bullet \ar@{}[l]|<<{10} \ar@{-}[lld]
  \\
  &&&&
  \bullet \ar@{}[r]|<<<{-5} \ar@{-}[rrd]&
  \bullet \ar@{}[r]|<<<{-3} \ar@{-}[rd]&
  &
  \bullet \ar@{}[l]|<<{3} \ar@{-}[ld]&
  \bullet \ar@{}[l]|<<{5} \ar@{-}[lld]
  \\
  &&&&&&
  \bullet \ar@{}[l]|<<{0}&
}
\]
 \caption{Spheres in the Cayley graph of $\Z$ with respect to $S=\{3,5\}$}
 \label{f:z-dead}
\end{figure}
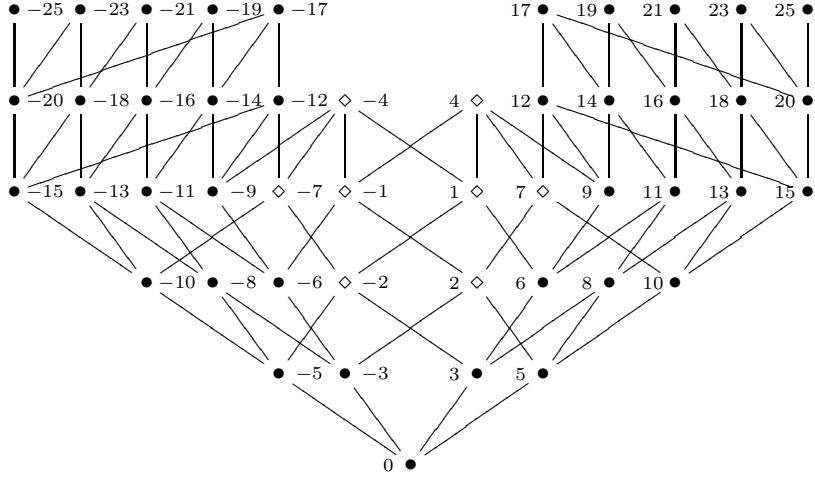

It is clear that $\pm 4$ are the only dead end elements in
$\Gamma(\Z,S)$. Note that the Frobenius values with respect to $S$
are $\pm 1,\pm 2,\pm 4, \pm 7$ (indicated by $\diamond$ in the
Cayley graph). While 7 is the largest Frobenius value, i.e., the
Frobenius number for $S$, it is clear that $\pm 4$ also play a
special role among Frobenius values.
\end{example}

The Cayley graph $\Gamma(G,S)$ induces a partial order on $G$.
Namely, $g \leq h$ if there exists a geodesic in $\Gamma(G,S)$ from
$1$ to $h$ that passes through $g$. Call this order on $G$ the
\emph{Cayley order} (with respect to $S$). Dead ends are precisely
the maximal elements in $G$ with respect to the Cayley order.

We are specifically interested in the relation between Frobenius
values and dead ends in $\Z$. Of course, when $S=\{1\}$ neither
Frobenius values nor dead ends exist. In the case when $S$ consists
of 2 elements we completely describe the connection.

\begin{theorem}\label{t:equiv}
Let $S=\{a,b\}$, where $a$ and $b$ are positive and relatively prime
integers. For an integer $d$, the following are equivalent:

(i) $d$ is a dead end with respect to $S$

(ii) $d$ is a maximal (in the Cayley order) Frobenius value with
respect to $S$.
\end{theorem}

The question of existence of dead ends in $\Z$ is touched upon
in~\cite{riley-w:notinvariant}, where Riley and Warshall show that
$m$ is a dead end with respect to both $\{2m,2m+1\}$ and
$\{2m-1,2m\}$ (with the exception of $m=1$ in the latter case). We
explicitly describe all dead ends in $\Z$, for any generating set
$S$ consisting of two elements.

\begin{theorem}\label{t:explicit}
Let $S=\{a,b\}$, where $a$ and $b$ are relatively prime integers
with $a > b \geq 1$. The dead ends in $\Z$ (the maximal Frobenius
values) with respect to $S$ are given as follows.

(i) If $a+b$ is even, then there are exactly $b-1$ dead ends, they
are all strict, they all have length $(a+b)/2$, and they are given
by
\[ d = \frac{(a+b)(2\alpha-b)}{2}, \]
for $\alpha=1,2,\dots,b-1$.

(ii) If $a+b$ is odd, then there are exactly $2(b-1)$ dead ends,
none of them is strict, they all have length $(a+b-1)/2$, and they
are given by
\[ d = \frac{(a+b)(2\alpha-b) \pm b}{2}, \]
for $\alpha=1,2,\dots,b-1$.
\end{theorem}

Note that Theorem~\ref{t:explicit} implies that there are no dead
ends when one of the generators is equal to $1$ (it is immediately
obvious that there are no Frobenius values in this case).

There are many examples of groups that have infinitely many dead
ends with respect to some generating sets. Such examples can be
found in~\cite{bogopolskii:bi-lipschitz} (the triangle group
$T(3,3,3)=\langle a, b | a^3=b^3=(ab)^3=1
\rangle$),~\cite{cleary-t:thompson-cayley} (Thompson's group
$F$),~\cite{cleary-t:lamplighter-dead}
(lamplighter),~\cite{cleary-r:fp-dead,riley-w:notinvariant} (some
finitely presented examples),~\cite{warshall:lattices} (some
lattices in $Nil$ and $Sol$, namely the discrete Heisenberg group
and any extension of $\Z^2$ by a hyperbolic automorphism). We show
that, on the contrary, $\Z$ can only have finitely many dead ends.

\begin{theorem}\label{t:finitelymany}
Let $S$ be a finite generating set of $\Z$. There exists only
finitely many dead ends in $\Z$ with respect to $S$.
\end{theorem}

Recall that, for a dead end $d$ of length $n$, the \emph{depth} of
$d$ is its distance to the sphere $\Sigma(n+1)$ decreased by 1 (some
authors prefer not to subtract 1 here, but this hardly matters).
Note that, in general, the condition that a group has only finitely
many dead ends is stronger than the condition that the depth of the
dead ends is uniformly bounded. For instance, the triangle group
$T(3,3,3)$ and Thompson's group $F$ have infinitely many dead ends
with respect to their standard 2-generator sets even though the
depth is uniformly bounded (see~\cite{bogopolskii:bi-lipschitz}
and~\cite{cleary-t:thompson-cayley}).

In the appendix we address the existence of dead ends in arbitrary
groups and show that every finitely generated group has a generating
set with respect to which dead ends exist (see
Theorem~\ref{t:theyexist} and Proposition~\ref{p:length2}).


\section{Proofs}

Let $a$ and $b$ be positive and relatively prime integers. Every
integer $c$ can be written as an integer liner combination $c=\alpha
a + \beta b$ of $a$ and $b$ in infinitely many ways. More precisely,
if $c=\alpha a + \beta b$ is one such representation then all other
integer solutions to the equation $c=xa+yb$ are given by
\begin{equation}\label{e:solutions}
 x = \alpha +bt, \qquad y = \beta - at,
\end{equation}
for integer values of $t$. There is a unique solution with $0 \leq x
< b$, which we call the $a$-normal form of $c$, and a unique
solution with $0 \leq y < a$, which we call the $b$-normal form of
$c$.

Observe that, for any finite generating set $S$ of $\Z$, the map $c
\mapsto -c$ is an automorphism of the Cayley graph of $\Z$ with
respect to $S$ and this automorphism fixes $0$. Therefore $c$ and
$-c$ have the same length with respect to $S$, and $c$ is a dead end
if and only if so is $-c$. Similarly, $c$ is a Frobenius value with
respect to $S$ if and only if so is $-c$, and $c$ is a maximal
Frobenius value if and only so is $-c$. We will freely use this
symmetry in the proofs that follow.

\subsection{Recognizing Frobenius values}

We begin by recalling a well known condition characterizing the
positively generated integers with respect to $S=\{a,b\}$ (it
appears, for instance, in~\cite{nijenhuis-w:representations}).

\begin{lemma}\label{l:who-is-pg}
Let $S=\{a,b\}$, where $a$ and $b$ are positive and relatively prime
integers, and let $c$ be a positive integer. The following
conditions are equivalent:

\textup{(i)} $c$ is positively generated with respect to $S$.

\textup{(ii)} the $a$-normal form $c = \alpha a + \beta b$ satisfies
the condition $\beta \geq 0$.
\end{lemma}

The characterization in Lemma~\ref{l:who-is-pg} will not be directly
useful to us (note that imposing the condition $\beta < 0$ in the
$a$-normal form lumps together all positive Frobenius values and all
negative integers), but the following slight modification will.

\begin{lemma}\label{l:who-is-f}
Let $S=\{a,b\}$, where $a$ and $b$ are positive and relatively prime
integers, and let $c$ be an integer (not necessarily positive). The
following conditions are equivalent:

\textup{(i)} $c$ is a Frobenius value with respect to $S$.

\textup{(ii)} the $a$-normal form $c = \alpha a + \beta b$ satisfies
the condition $-a < \beta < 0 < \alpha < b$.

\textup{(iii)} the $b$-normal form $c = \alpha a + \beta b$
satisfies the condition $-b < \alpha < 0 < \beta < a$.
\end{lemma}

\begin{proof}
(i) implies (ii). Let $c$ be a Frobenius value with respect to $S$
and let $c = \alpha a + \beta b$ be the $a$-normal form of $c$.
Since $c$ is a Frobenius value neither $\alpha$ nor $\beta$ can be
equal to $0$. Thus $0< \alpha < b$ and $\beta$ must be negative. If
$\beta \leq -a$ then $c = (\alpha-b)a + (\beta +a)b$ and since
$\alpha - b <0$ and $\beta+a \leq 0$ we obtain that $c$ is
negatively generated, a contradiction. Thus $-a < \beta < 0$.

(ii) implies (i). Let the $a$-normal form $c = \alpha a + \beta b$
satisfy the condition $-a < \beta < 0 < \alpha < b$. All other
solutions to the equation $c = xa+yb$ are given
by~\eqref{e:solutions}. For $t \geq 0$, $x=\alpha + bt \geq \alpha
>0$ and $y = \beta - at \leq \beta <0$, while for $t \leq -1$,
$x=\alpha + bt \leq \alpha -b  < 0$ and $y = \beta - at \geq \beta
+a >0$. Thus, in any representation of $c$ in the form $c=xa+yb$,
one of the integers $x$ and $y$ is positive while the other is
negative, implying that $c$ is a Frobenius value with respect to
$S$.

(ii) is equivalent to (iii). If $c = \alpha a + \beta b$ and $-a <
\beta < 0 < \alpha < b$, then $c = (\alpha - b)a + (\beta + a)b$ and
$-b < \alpha - b < 0 < \beta + a < a$. Thus (ii) implies (iii) and,
by symmetry, (iii) implies (ii).
\end{proof}

It is interesting to observe that Lemma~\ref{l:who-is-f} provides a
rather simple proof of the following classical result of
Sylvester~\cite{sylvester:number-ab}.

\begin{corollary}
Let $S=\{a,b\}$, where $a$ and $b$ are positive and relatively prime
integers. The number of positive Frobenius values with respect to
$S$ is equal to $(a-1)(b-1)/2$.
\end{corollary}

\begin{proof}
Since Frobenius values are exactly the numbers $\alpha a + \beta b$,
for $-a < \beta < 0 < \alpha < b$, and no two such numbers are equal
(every integer has a unique $a$-normal form) the number of Frobenius
values is $(a-1)(b-1)$. By symmetry, exactly half of them are
positive.
\end{proof}

Three proofs of Sylvesters's result are offered
in~\cite{ramirez-alfonsin:frobenius-book}, but none of them uses the
above argument for the simple reason that they all only consider and
concentrate on positive Frobenius values. For instance, one of the
proofs offered in~\cite{ramirez-alfonsin:frobenius-book} follows the
argument of Nijenhuis and Wilf~\cite{nijenhuis-w:representations}
(and is based on Lemma~\ref{l:who-is-pg}) and shows that an integer
$c$ in the closed interval $[0,ab-a-b]$ is a Frobenius value if and
only if $ab-a-b-c$ is not.


\subsection{Recognizing maximal Frobenius values}

We now concentrate on description of all maximal Frobenius values.

\begin{lemma}\label{l:length}
Let $S=\{a,b\}$, where $a$ and $b$ are positive and relatively prime
integers. Let $c = \alpha a + \beta b$ be the $a$-normal form of the
Frobenius value $c$. The length of $c$ with respect to $S$ is
achieved either at the $a$-normal form or at the $b$-normal form of
$c$, i.e.,
\[ \ell(c) = \min \{ |\alpha|+|\beta|, a+b-(|\alpha|+|\beta|)\}. \]
\end{lemma}

\begin{proof}
Consider the solutions~\eqref{e:solutions} to the equation $c =
xa+yb$. We have
\[
 |x|+|y| = |\alpha + bt| +|\beta - at| =
 \begin{cases}
 (a+b)t + |\alpha|+|\beta|,     & t \geq 0 \\
 (a+b)|t| - (|\alpha|+|\beta|), & t \leq -1
 \end{cases}.
\]
Thus the length of $c$ is achieved at $t=0$, which is the $a$-normal
form of $c$, or at $t=-1$, which is the $b$-normal form of $c$, and
$\ell(c) = \min \{ |\alpha|+|\beta|, a+b-(|\alpha|+|\beta|)\}$.
\end{proof}

Note that Lemma~\ref{l:length} implies that there are no Frobenius
values of length greater than $\left\lfloor \frac{a+b}{2}
\right\rfloor$. We describe explicitly the Frobenius values of
length exactly $\left\lfloor \frac{a+b}{2} \right\rfloor$.

\begin{lemma}\label{l:explicit}
Let $S=\{a,b\}$, where $a$ and $b$ are relatively prime integers
with $a > b \geq 1$. The Frobenius values of length $\left\lfloor
\frac{a+b}{2} \right\rfloor$ with respect to $S$ are given as
follows.

\textup{(i)} If $a+b$ is even, then there are exactly $b-1$
Frobenius values of length $(a+b)/2$ and they are given by
\[ c = \alpha a + \left( \alpha -  \frac{(a+b)}{2} \right) b, \]
for $\alpha=1,2,\dots,b-1$.

\textup{(ii)} If $a+b$ is odd, then there are exactly $2(b-1)$
Frobenius values of length $(a+b-1)/2$ and they are given by
\[ c = \alpha a + \left( \alpha -  \frac{(a+b \mp 1)}{2} \right) b, \]
for $\alpha=1,2,\dots,b-1$.
\end{lemma}

\begin{proof}
The claim easily follows from Lemma~\ref{l:who-is-f} and
Lemma~\ref{l:length}.

(i) Note that, when $0< \alpha < b$,
\[
 -a < - \frac{a+b}{2} < \alpha - \frac{a+b}{2} < b - \frac{a+b}{2} =
 \frac{b-a}{2} < 0.
\]
Thus $\beta=\alpha - (a+b)/2$ satisfies the condition $-a < \beta <
0$ and therefore $c = \alpha a + \beta b$ is a Frobenius value. The
length of $c$ is $|\alpha|+|\beta| = \alpha-\beta = (a+b)/2$.

Conversely, if $c = \alpha a + \beta b$ is a Frobenius value of
length $(a+b)/2$ and $0<\alpha$ then $\beta$ must be negative and we
must have $(a+b)/2=|\alpha|+|\beta| = \alpha - \beta$, which implies
$\beta = \alpha - (a+b)/2$.

(ii) The proof is analogous to the one given for (i). The difference
in the number of solutions comes from the fact that if $|\alpha| +
|\beta|$ is equal to either $(a+b-1)/2$ or $(a+b+1)/2=
a+b-(a+b-1)/2$, then the length of the corresponding Frobenius value
$c= \alpha a + \beta b$ is $(a+b-1)/2$.
\end{proof}

The next result shows that Lemma~\ref{l:explicit} describes exactly
the maximal Frobenius values.

\begin{lemma}\label{l:equiv}
Let $S=\{a,b\}$, where $a$ and $b$ are positive and relatively prime
integers. A Frobenius value $c$ with respect to $S$ is maximal
Frobenius value if and only if its length is $\left\lfloor
\frac{a+b}{2} \right\rfloor$.
\end{lemma}

\begin{proof}
If a Frobenius value $c$ has length $\left\lfloor \frac{a+b}{2}
\right\rfloor$ then it is certainly maximal Frobenius value, since
the length of a maximal Frobenius value cannot be greater than
$\left\lfloor \frac{a+b}{2} \right\rfloor$ (by
Lemma~\ref{l:length}).

Conversely, let $c$ be a Frobenius value of length strictly smaller
than $\left\lfloor \frac{a+b}{2} \right\rfloor$ and let $c = \alpha
a + \beta b$ be its $a$-normal form.

Let the length of $c$ be achieved at the $a$-normal form, i.e., let
$\ell(c) = |\alpha| + |\beta|$. We cannot have $\alpha=b-1$ and
$\beta=1-a$ (otherwise $\ell(c) = a+b-2 \geq (a+b)/2 +3/2-2 =
(a+b-1)/2$, a contradiction). Thus either $c'=c+a$ (if $\alpha \neq
b-1$) or $c'=c-b$ (if $\beta \neq 1-a$) is a Frobenius value. Since
$|\alpha|+|\beta|+1 \leq \left\lfloor \frac{a+b}{2} \right\rfloor$
we in each case obtain that $\ell(c')= |\alpha|+|\beta|+1 =
\ell(c)+1$. Thus $c$ is not a maximal Frobenius value.

If the length of $c$ is achieved at the $b$-normal form then the
length of $-c$ is achieved at its $a$-normal form (multiplying the
$b$-normal form of $c$ by $-1$ throughout provides the $a$-normal
form of $-c$). By the previous argument, $-c$ is not a maximal
Frobenius value and, by symmetry, $c$ is not a maximal Frobenius
value either.
\end{proof}


\subsection{Connection to dead ends}

We now provide the proofs of the statements relating the maximal
Frobenius values and the dead ends. In fact, we provide two proofs
that maximal Frobenius values are dead ends. One is based on the
fact that we already know explicitly all maximal Frobenius values
and can relatively easily and directly check that the length of any
of their neighbors in the Cayley graph does not exceed $\left\lfloor
\frac{a+b}{2} \right\rfloor$, which is the length of the maximal
Frobenius values. The other proof is more conceptual and relies
solely on the definition of a maximal Frobenius value (it does not
use Lemma~\ref{l:who-is-f}, Lemma~\ref{l:explicit}, or
Lemma~\ref{l:equiv}). While the first proof is slightly shorter and
contributes to the proof of Theorem~\ref{t:explicit}, the second
proof is more likely to be amenable to generalizations to larger
generating sets (since it is unlikely that explicit descriptions of
maximal Frobenius values for such generating sets can be found).

\begin{proof}[First proof of Theorem~\ref{t:equiv}, \textup{(ii)} implies \textup{(i)}]
Let $c$ be a maximal Frobenius value. Without loss of generality,
assume that the length of $c$ is achieved at the $a$-normal form
$c=\alpha a+ \beta b$ (otherwise we may consider $-c$). Thus $-a <
\beta < 0 < \alpha < a$ and $\ell(c) = |\alpha|+ |\beta| =
\left\lfloor \frac{a+b}{2} \right\rfloor$.

Since
\[
 c-a = (\alpha-1)a + \beta b \qquad \text{and}\qquad c+b = \alpha a
+ (\beta+1) b
\]
and $|\alpha -1|+|\beta| = |\alpha|+|\beta|-1 = |\alpha|+|\beta+1|$,
we see that $\ell(c-a),\ell(c+b) \leq \ell(c)-1$ (and therefore
$\ell(c-a)= \ell(c+b) = \ell(c)-1$).

Further,
\begin{align*}
 c+a &= (\alpha+1)a + \beta b = (\alpha+1-b)a + (\beta+a) b, \\
 c-b &= \alpha a + (\beta -1) b = (\alpha-b)a + (\beta-1+a) b,
\end{align*}
and $|\alpha+1-b| + |\beta+a| = a+b-(|\alpha|+|\beta|)-1 =
|\alpha-b| + |\beta-1+a|$. Therefore, if $a+b$ is even, $\ell(c+a),
\ell(c-b) \leq a+b-(|\alpha|+|\beta|)-1  = \frac{a+b}{2} -1 =
\ell(c)-1$ (implying that $\ell(c+a), \ell(c-b)=\ell(c)-1$) and, if
$a+b$ is odd, $\ell(c+a), \ell(c-b) \leq a+b-(|\alpha|+|\beta|)-1 =
\frac{a+b+1}{2} -1 = \ell(c)$.
\end{proof}

\begin{proof}[Second proof of Theorem~\ref{t:equiv}, \textup{(ii)} implies \textup{(i)}]
Assume that $c$ is a maximal Frobenius value. Without loss of
generality, assume that $c=\alpha a+ \beta b$, with $\beta < 0 <
\alpha$, and $\ell(c)=|\alpha|+|\beta|$ (note that neither $\alpha$
nor $\beta$ can be 0, since $m$ is a Frobenius value).

Since
\[
 c-a = (\alpha-1)a + \beta b \qquad \text{and}\qquad c+b = \alpha a
+ (\beta+1) b
\]
and $|\alpha -1|+|\beta| = |\alpha|+|\beta|-1 = |\alpha|+|\beta+1|$,
we see that $\ell(c-a),\ell(c+b) \leq \ell(c)-1$ (and therefore
$\ell(c-a)= \ell(c+b) = \ell(c)-1$).

Consider $c-b$ and $c+a$. Let $c-b = \alpha'a+\beta'b$, with
$\ell(c-b)=|\alpha'|+|\beta'|$. Then $c=\alpha'a+(\beta'+1)b$,
$c+a=(\alpha'+1)a+(\beta'+1)b$, and since $c$ is a Frobenius value,
either $\beta'+1<0<\alpha'$ or $\alpha'<0<\beta'+1$. In the former
case
\[
 \ell(c+a) \leq |\alpha'+1|+|\beta'+1| = (\alpha'+1)-(\beta'+1) =
  \alpha'-\beta' = |\alpha'|+|\beta'| = \ell(c-b)
\]
and in the latter
\[
 \ell(c+a) \leq |\alpha'+1|+|\beta'+1| = -(\alpha'+1)+(\beta'+1) =
  -\alpha'+\beta' = |\alpha'|+|\beta'| = \ell(c-b).
\]
Thus $\ell(c+a) \leq \ell(c-b)$ and, by symmetry, it follows that
$\ell(c-b)=\ell(c+a)$.

Assume that $c$ is not a dead end. Then both $c-b$ and $c+a$ have
length $\ell(c)+1$. Since $c$ is a maximal Frobenius value, none of
$c-b$ and $c+a$ is a Frobenius value. However, $c+a$ cannot be
negatively generated since $c = (c+a)-a$ would then also be
negatively generated, which contradicts the assumption that $c$ is a
Frobenius value. If $c+a = \alpha''a + \beta''b$, where
$\alpha'',\beta''\geq 0$, then $\alpha''=0$ (if $\alpha''>0$ then
$c=(\alpha''-1)a+\beta'' b$ would be positively generated). Thus
$c+a = \beta''b$, for some positive $\beta''$. Moreover, since $c$
is a Frobenius value, its length is at least 2, implying $\ell(c+a)
\geq 3$ and therefore $\beta'' \geq 3$. By a symmetric argument,
$c-b$ must be negatively generated, and $c-b=-\alpha''a$, for some
$\alpha'' \geq 3$. But then
\[ a+b = (c+a)-(c-b) = \beta''b + \alpha''a \geq 3b + 3a = 3(a+b), \]
a contradiction.
\end{proof}

\begin{proof}[Proof of Theorem~\ref{t:equiv}, \textup{(i)} implies \textup{(ii)}]
Without loss of generality, assume that $a > b$.

Let $c$ be a positive integer and let $c=\alpha a + \beta b$, where
$\alpha$ and $\beta$ are integers (we are not assuming any normal
form here). All other solutions to the equation $c=xa+yb$ are given
by~\eqref{e:solutions}. The length of $c$ with respect to $S$ is the
minimal value of the function
\[ f_c(t) = |\alpha+bt| + |\beta-at|, \]
at an integer value of $t$. We have
\[ f_c(t) = \begin{cases}
            -\alpha+\beta-(a+b)t, & t \leq -\alpha/b, \\
            c/b, & t=-\alpha/b, \\
            \alpha+\beta-(a-b)t, & -\alpha/b \leq t \leq \beta/a,\\
            c/a, & t = \beta/a, \\
            \alpha-\beta+(a+b)t, & \beta/a \leq t
            \end{cases}
\]
and the graph of $f_c(t)$ is given in Figure~\ref{ca-weight} (full
line).
\begin{figure}[!h]
\begin{center}
 \includegraphics[height=170pt]{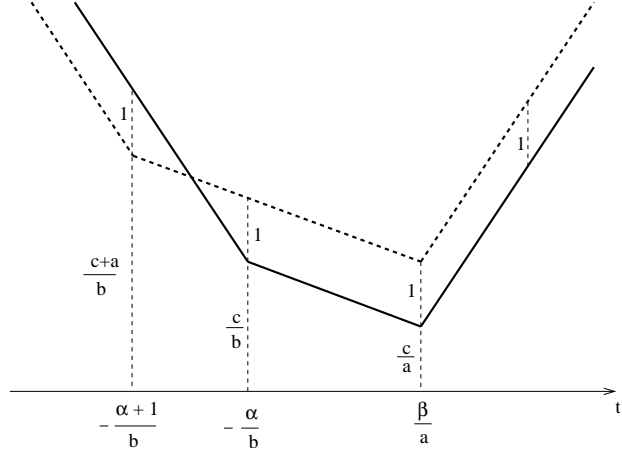}
\end{center}
\caption{Graphs of the functions $f_c(t)$ and
$f_{c+a}(t)$}\label{ca-weight}
\end{figure}

Since the function $f_c(t)$ is decreasing for values of $t$ smaller
than $\beta/a$, achieves its minimum at $\beta/a$ and is increasing
for values of $t$ greater than $\beta/a$, the length of $c$ is
obtained either at the integer $\ell=\left\lfloor \beta/a
\right\rfloor$ that is closest to $\beta/a$ to the left of $\beta/a$
or the integer $r=\lceil \beta/a \rceil$ that is closest to
$\beta/a$ to the right of $\beta/a$.

Consider the function $f_{c+a}(t)$ that determines the length of
$c+a$. We have $c+a = (\alpha+1)a + \beta b$ and the graph of the
function $f_{c+a}(t)$ is given in Figure~\ref{ca-weight} (dotted
line). If there exists an integer in the closed interval
$[-\alpha/b, \beta/a]$, then $c$ cannot be a dead end (to the right
of $-\alpha/b$ the function $f_{c+a}(t)$ is 1 unit above $f_c(t)$,
implying that the length of $c+a$ is larger than the length of $c$).

Assume that $d$ is a positive dead end, $d=\alpha a+\beta b$, and
$\alpha$ and $\beta$ are chosen so that $\ell(d)=|\alpha|+|\beta|$.
This means that the minimum of $f_d(t)$ at an integer point is
achieved at $t=0$ and there is no integer in the interval
$[-\alpha/b, \beta/a]$. Thus, either $-1 < -\alpha/b < \beta/a < 0$
or $0 < -\alpha/b < \beta/a < 1$. In the former case $-a < \beta < 0
< \alpha < b$, while in the latter $-b < \alpha < 0 < \beta < a$. In
each case we conclude that $d$ is a Frobenius value by
Lemma~\ref{l:who-is-f}.

If $d$ is a negative dead end, then $d$ is just negative of some
positive dead end, so $d$ is a Frobenius value in this case as well.

Thus we proved that all dead ends are Frobenius values. Since dead
ends are maximal in the Cayley order they must be maximal Frobenius
values as well.
\end{proof}

\begin{proof}[Proof of Theorem~\ref{t:explicit}]
By Lema~\ref{l:equiv}, the maximal Frobenius values are precisely
the Frobenius values of length $\left\lfloor \frac{a+b}{2}
\right\rfloor$, and these values are explicitly described in
Lemma~\ref{l:explicit}. By Theorem~\ref{t:equiv} these values are
precisely the dead ends in $\Z$ with respect to $S$.

Since
\[
 \alpha a + \left( \alpha - \frac{(a+b)}{2} \right) b =
     \frac{(a+b)(2\alpha-b)}{2},
\]
and
\[
 \alpha a + \left( \alpha - \frac{(a+b \mp 1)}{2} \right) b =
     \frac{(a+b)(2\alpha-b) \pm b}{2},
\]
for $\alpha=1,2,\dots,b-1$, the lists of numbers given in the claims
of Lemma~\ref{l:explicit} and Theorem~\ref{t:explicit} coincide.

The claim on the lengths of dead ends follows from the corresponding
claim on the lengths of maximal Frobenius values
(Lemma~\ref{l:equiv}).

In the course of the first proof of the direction (ii) implies (i)
of Theorem~\ref{t:equiv} we already proved that, if $a+b$ is even
then, for any dead end $c$, $\ell(c \pm a) = \ell(c \pm b) =
\ell(c)-1$. Thus, in this case all dead ends are strict.

The equality
\[
 \frac{(a+b)(2\alpha-b) + b}{2} - b =
 \frac{(a+b)(2\alpha-b) - b}{2},
\]
for $\alpha = 1,\dots,b-1$, shows that, when $a+b$ is odd, each dead
end has another dead end as a neighbor in the Cayley graph, which
then means that none of them is a strict dead end (since they all
have the same length).
\end{proof}

\begin{proof}[Proof of Theorem~\ref{t:finitelymany}]
Without loss of generality, let $S=\{b_1,b_2,\dots,b_k\}$ be a
generating set for $\Z$ such that $k \geq 1$ and $b_1 > b_2 > \dots
> b_k >0$. Denote $b_1=a$.

Let
\[ c = x_1a + x_2b_2 + \dots + x_kb_k, \]
and let the length $\ell(c)$ of $c$ with respect to $S$ be given by
\[ \ell(c) = |x_1| + |x_2| + \dots + |x_k|. \]
We claim that $|x_i|<a$, for $i=2,\dots,k$. Indeed, assume $x_i \geq
a$, for some $i=2,\dots,k$. Then we have
\[ c = (x_1+b_i)a + x_2b_2 + \dots + (x_i-a)b_i + \dots + x_kb_k,\]
which implies that
\[ \ell(c) \leq |x_1+b_i|+|x_2|+ \dots + |x_i-a| + \dots + |x_k| \leq |x_1|+\dots+|x_k| + b_i - a < \ell(c),\]
a contradiction. Thus we must have $x_i < a$, for $i=2,\dots,k$. A
symmetric argument shows that $-a<x_i$, for $i=2,\dots,k$.

We show that if $c > (a-1)b$, where $b=b_2 + \dots + b_k$, then
$x_1$ must be positive. Indeed, if $x_1$ is not positive then
\[ c = x_1a + x_2b_2 + \dots + x_kb_k \leq x_2b_2 + \dots + x_kb_2 \leq (a-1)(b_2+ \dots + b_k) = (a-1)b.\]
Thus, by contraposition, if $c$ is large enough $x_1$ must be
positive.

We show that if $c>(a-1)b$ then $c$ cannot be a dead end. Indeed, in
that case $c+a$ is also large enough, so that if
$c+a=x_1'a+x_2'b_2+\dots+ x_k'b_k$, with $\ell(c) = |x_1'| + |x_2'|
+ \dots + |x_k'|$, then $x_1'$ is positive. But then
$c=(x_1'-1)a+x_2'b_2+\dots +x_k'b_k$ and therefore
\[ \ell(c)\leq |x_1'-1|+|x_2'| +\dots + |x_k'| = |x_1'|+|x_2'| +\dots + |x_k'| -1  = \ell(c+a)-1,\]
which shows that $c$ is not a dead end.

Thus there are only finitely many positive dead ends and, by
symmetry, there are only finitely many dead ends in $\Z$.
\end{proof}

\section*{Appendix: all groups have dead ends}

Warshall showed~\cite{warshall:many} that if $G$ is a finitely
generated group with infinitely many finite homomorphic images then,
for every $n$, $G$ has a finite generating set with respect to which
$G$ has dead ends of depth at least $n$.

We show that, surprisingly, if one is only interested in existence
of dead ends, no conditions on the group are needed (for
applications involving just the existence of dead ends
see~\cite{bartholdi:amenablealg2}).

\begin{theorem}\label{t:theyexist}
Every infinite, finitely generated group $G$ has a finite generating
set $S$ with respect to which $G$ has dead ends.

If $G$ is generated by $m$ elements, the size of $S$ can be chosen
to be no greater than $4m+2$.
\end{theorem}

\begin{proof}
Let $G$ have an element $a$ of order at least 5 (including the
possibility of infinite order). Let $S''$ be any generating set for
$G$ and let $S'$ be the set obtained from $S''$ by removing any
generators that happen to be in $\langle a \rangle$. For every $b$
in $S'$ define
\[ S_b = \{b,ab,ab^{-1},aba^{-1}\}, \]
and let
\[ S= \left( \bigcup_{b \in S'} S_b \right) \cup \{a^2, a^3 \}. \]

It is clear that $S$ is a generating set for $G$ (of size at most
$4|S''|+2$). We claim that $a$ is a dead end of length 2 with
respect to $S$.

Since $a = a^3a^{-2}$ the length of $a$ is no greater than 2. On the
other hand, $a$ is not equal to any of the generators in $S$ or
their inverses (this is because the order of $a$ is at least 5 and
$b \not \in \langle a \rangle$). Thus $a$ has length 2.

We have
\begin{alignat*}{2}
 a \cdot a^2 & = a^3 & a \cdot a^{-2} & = a^2 \cdot a^{-3} \\
 a \cdot a^3 & = a^2 \cdot a^2 & a \cdot a^{-3} &= a^{-2} \\
 a \cdot b & =  ab & a \cdot b^{-1} & = ab^{-1} \\
 a \cdot ab &  = a^2 \cdot b & a \cdot (ab)^{-1} & = (aba^{-1})^{-1} \\
 a \cdot ab^{-1} &  =  a^2 \cdot b^{-1} & a \cdot (ab^{-1})^{-1} & = aba^{-1} \\
 a \cdot aba^{-1} &  =  a^2 \cdot (ab^{-1})^{-1} & \qquad\qquad a \cdot (aba^{-1})^{-1} &  = a^2 \cdot (ab)^{-1},
\end{alignat*}
which shows that no neighbor of $a$ in the Cayley graph of $G$ with
respect to $S$ has length higher than 2. Therefore $a$ is a dead end
of length 2.

Let $G$ have an element $a$ of order 4. Define $S'$ as before. A
similar argument to the one given above then shows that $a^2$ is a
dead end of length 2 with respect to the generating set $S$ defined
by
\[ S= \left( \bigcup_{b \in S'} S_b \right) \cup \{a \}. \]

Finally, if $G$ has no elements of order $4$ or higher, then $G$ is
a finitely generated group in which $g^6=1$, for all $g \in G$.
Therefore $G$ is a homomorphic image of the free Burnside group
$B(m,6)$. Since $B(m,6)$ is finite~\cite{mhall:bm6}, this means that
$G$ is finite, contradicting the assumption that $G$ is infinite.
\end{proof}

The statement in Theorem~\ref{t:theyexist} is concerned only with
infinite groups since finite groups always have dead ends. However,
for completeness, we observe that dead ends of length 2 can be
achieved in any group that is sufficiently large to allow elements
of length 2 to exist (thus all groups but the trivial group and the
cyclic groups of order 2 and 3).

\begin{proposition}\label{p:length2}
Every finitely generated group $G$ that has at least 4 elements has
a finite generating set $S$ with respect to which $G$ has a dead end
of length 2.
\end{proposition}

\begin{proof}
The proof of Theorem~\ref{t:theyexist} applies unless $G$ is a
finite group in which every element has order 2 or 3. In that case,
let $a$ be any nontrivial element of $G$. Since $a$ has order 2 or 3
and $|G|\geq 4$ the set $\{1,a,a^{-1}\}$ is a proper subgroup of
$G$. The set $S=G-\{1,a,a^{-1}\}$ is then a finite generating set
for $G$ (since the complement of a proper subgroup always generates
the group). The element $a$ is then a dead end of length 2 with
respect to $S$.
\end{proof}

\def\cprime{$'$}

\end{document}